\documentclass[10pt]{article}
\usepackage{amssymb, amsmath, amsthm}

\newtheorem{theorem}{Theorem}[section]
\newtheorem{lemma}[theorem]{Lemma}
\newtheorem{proposition}[theorem]{Proposition}
\newtheorem{corollary}[theorem]{Corollary}

\newtheorem{definition}[theorem]{Definition}

\newtheorem{remark}[theorem]{Remark}

\numberwithin{equation}{section}

\newcommand\NN{\mathbb{{N}}}

\newcommand\RR{\mathbb{{R}}}
\newcommand\CC{\mathbb{{C}}}
\newcommand\ZZ{\mathbb{{Z}}}

\newcommand\rE{{\rm E}}

\newcommand\ra{{\rm a}}

\newcommand\rd{{\rm d}}

\newcommand\rp{{\rm p}}
\newcommand\rrq{{\rm q}}

\newcommand\be{{\boldsymbol e}}

\newcommand\bh{{\bf h}}
\newcommand\bj{{\bf j}}

\newcommand\bell{{\boldsymbol\ell}}

\newcommand\bt{{\bf t}}

\newcommand\bx{{\boldsymbol x}}

\newcommand\bz{{\boldsymbol z}}

\newcommand\balpha{{\boldsymbol\alpha}}
\newcommand\bbeta{{\boldsymbol\beta}}
\newcommand\bdelta{{\bf \delta}}
\newcommand\bveps{{\boldsymbol\epsilon}}
\newcommand\beps{{\boldsymbol\varepsilon}}

\newcommand\btau{{\boldsymbol\tau}}
\newcommand\bxi{{\boldsymbol\xi}}

\newcommand\bnull{{\boldsymbol0}}
\newcommand\bone{{\boldsymbol1}}

\newcommand\cS{{\mathcal S}}

\newcommand{\mfrac}[2]%
{\raisebox{0.5pt}{\footnotesize$\dfrac{#1}{#2}$}}
\newcommand{\mbinom}[2]%
{\raisebox{0.5pt}{\footnotesize$\dbinom{#1}{#2}$}}
\def\smmat\{#1&#2\cr#3&#4\}%
   {\!\!\pmatrix{\!#1\!&\!#2\!\cr\!#3\!&\!#4\!}\!\!}
\newcommand\scrm{{\raise0.5pt\hbox{-}}}
\def\eop{{ \vrule height7pt width7pt depth0pt}\par\bigskip} 
\newcommand\ie{{\it\thinspace i.e.}}
\newenvironment{pf}%
{\par\noindent\textbf{Proof:}~}

\begin{document}

\begin{titlepage}
\title{Polynomial Reproduction of Multivariate Scalar Subdivision Schemes}
\author{Maria Charina and Costanza Conti}
\end{titlepage}

\maketitle
\begin{abstract}
A stationary subdivision scheme generates the full space of polynomials of degree up to $k$ if and
only if its mask satisfies sum rules of order $k+1$, or its symbol satisfies zero conditions
of order $k+1$. This property is often called the polynomial reproduction property of the subdivision scheme.
It is a well-known fact that this property is, in general, only necessary for
the associated refinable function to have approximation order $k+1$.\\
In this paper we study a different polynomial reproduction property of
multivariate scalar subdivision scheme with dilation matrix $mI,\
|m| \ge 2$. Namely, we are interested in capability of a subdivision scheme to reproduce
in the limit exactly the same polynomials from which the data is sampled.
The motivation for this paper are the results in \cite{ALevin2003} that
state that such a reproduction property of degree $k$ of the subdivision scheme is sufficient for having
approximation order $k+1$. \\
Our main result yields simple
algebraic conditions on the subdivision symbol for computing the
exact degree of such polynomial reproduction and also for determining the associated parametrization. The
parametrization determines the grid points to which the newly
computed values are attached at each subdivision iteration to ensure the higher degree
of polynomial reproduction. We illustrate our results with several examples.

\end{abstract}


\noindent Keyword:
      Subdivision schemes, polynomial reproduction, subdivision parametrization, approximation order

\maketitle
\section*{Introduction}

Interest in subdivision schemes is motivated by their applications in computer graphics, computer
aided geometric design, animation, wavelet and frame construction.
Important properties of subdivision schemes such as convergence,
regularity, polynomial generation, approximation order, etc.,  have been
studied by several authors, see surveys \cite{Cabrelli2, CDM91, D92,
DynLevin02, JiaJiang, JiaJiang2} and references therein. In this paper we would like to
distinguish between the concepts of polynomial generation and polynomial reproduction
of subdivision schemes. The so-called polynomial generation of degree $k$ is the capability of subdivision
to generate the full space of polynomials of degree up to $k$. This property
is equivalent to sum rules of order $k+1$ on the subdivision mask, or, equivalently, to zero conditions of order $k+1$ on the subdivision symbol,
see e.g. \cite{CDM91, JePlo, JiaJiang}. Polynomial
generation of degree $k$ also implies that the associated refinable function has accuracy of order
$k$ \cite{JiaJiang2}, but is, in general, only necessary for the corresponding shift-invariant space to have
approximation order $k+1$. This has been already observed for B-splines and box splines \cite{BHR}. The so-called polynomial reproduction
is the capability of subdivision schemes
to produce in the limit exactly the same polynomials from which the data is sampled. The results in \cite[Section 2.4]{ALevin2003}
state that polynomial reproduction of degree $k$ of convergent subdivision is sufficient for the
associated shift-invariant space to have approximation order $k+1$. This motivates our interest in
polynomial reproduction of subdivision schemes.

Our main goal is to derive simple algebraic conditions on subdivision symbol that
allow us to determine the degree of its polynomial reproduction. Note that the concepts of polynomial reproduction
and generation coincide in the case $k=0$, \ie\ in the case of reproduction of constants. In the $L_2$ setting, polynomial generation of degree $k$ is also sufficient for
approximation order $k+1$, see e.g. \cite{JePlo}. For convergent interpolatory subdivision schemes, the concepts of
polynomial generation and reproduction are equivalent and, thus, characterize the approximation power of the
corresponding shift-invariant space, see \cite{Je}. We emphasize that
there is a multitude of results on approximation order of a refinable function.  Those results however are mostly
derived from the properties of the associated shift-invariant space and not from the
properties of the coefficients of the refinement equation -  the subdivision mask, in the subdivision context.  There are also data
pre-processing techniques for achieving the optimal approximation
order of a shift-invariant space associate with a convergent subdivision
scheme, \cite{ALevin2003, DynLevin02, Jia98}. The polynomial
reproduction of order $k$  makes preprocessing unnecessary, which is
undoubtedly an advantage.

In the univariate case, polynomial reproduction has been studied
in \cite{DynShen2008} for binary primal and dual subdivision schemes and
extended in \cite{ContiHormann} to univariate subdivision schemes
of any a-rity. In \cite{ContiHormann} the authors provide unified
algebraic conditions on the subdivision symbol for polynomial
reproduction with no restrictions on the associated
parametrization - the grid points to which the newly computed
values are attached at each subdivision iteration. The results of
our paper extend \cite{ContiHormann} to the multivariate setting
for scalar subdivision with dilation matrix $mI,\ |m| \ge 2$. To the
best of our knowledge our results are the first ones on multivariate
polynomial reproduction of subdivision schemes.  Our interest in the case of dilation
matrix $mI$ is motivated, e.g., by the bivariate $\sqrt{3}-$subdivision whose refinable
function is also refinable with respect to dilation $-3I$ and iterated
mask $a(z_1 z_2^{-2},z_1^2 z_2^{-1}) \cdot a(z_1,z_2)$. There are several other
expansive dilation matrices $M$ which satisfy $M^n=mI$, and, thus such that our results
are applicable.

\smallskip
The main result of our paper, Theorem \ref{teo:main},  states that for a non-singular
subdivision scheme with finitely supported mask $\ra = \{\ra_\balpha,\ \balpha\in\ZZ^s\}$
and symbol $\displaystyle{a(\bz)=\sum_{\balpha\in\ZZ^s} \ra_\balpha\;\bz^\balpha}$ the
polynomial reproduction of order $k$ is equivalent to
\begin{equation}\label{eq:main}
 \bigl(D^\bj a\bigr)(\bone)=|m|^s \prod_{i=1}^s \prod_{\ell_i=0}^{j_i-1}(\tau_i-\ell_i)
 \quad \hbox{and} \quad \bigl(D^\bj a\bigr)(\beps)=0 \ \
\hbox{for} \ \ \beps\in \Xi',\ \ |\bj| \le k\,,
\end{equation}
where $\Xi'$ is a finite set of certain multi-indices and
$\btau=(\tau_1,\cdots,\tau_s)\in \RR^s$ appears in the
parametrization associated with the subdivision scheme. The
importance of condition \eqref{eq:main} for $k=1$ is that it
allows us to identify the correct parametrization that for any
non-singular or for even only convergent subdivision scheme
guarantees at least the reproduction of linear polynomial. The parametrization
determines the grid points to which the newly computed values are attached at
each step of subdivision recursion to ensure the higher degree of
polynomial reproduction of a scheme.

\smallskip This paper is organized as follows. The first section sets the
notation, provides the background on multivariate subdivision
schemes stressing the difference between polynomial reproduction
and polynomial generation. We also define
sequence of parameter values associated with a subdivision scheme,
\ie, the subdivision parametrization. In section 2, we first
provide algebraic tools for determining the correct
parametrization needed to ensure reproduction of linear polynomials. We also give
the necessary and sufficient conditions on the symbol of any non-singular subdivision scheme
that guarantee polynomial reproduction
up to a certain degree. Effect of a shift of the mask  on the
degree of polynomial reproduction is investigated in Section 3.
There we also show that the concept of polynomial reproduction and
polynomial generation are equivalent for convergent
interpolatory schemes. Thus, reproducing the results in \cite{Je}. In section 3 we also provide the correct
parametrization for box spline subdivision schemes together with
several examples. The effect of the shifts of the box splines on the
approximation order of the corresponding shift-invariant spaces has been
already observed e.g. in \cite{BHR}.

\section{Background  and notation} 
\subsection{Subdivision scheme}
A scalar $s$-variate subdivision scheme with a \emph{dilation matrix} $mI$, $|m| \ge 2$, is given by a scalar
finitely supported sequence $\ra = \{\ra_\balpha\in \RR,\ \balpha\in\ZZ^s\}$,
the so-called \emph{mask}. The \emph{subdivision operator} $\cS_\ra$ acting on
data sequences $\rd=\{\rd_\balpha,\ \balpha\in\ZZ^s\} \in
\ell(\ZZ^s)$ is defined by
\begin{equation}\label{eq:subdop}
  \big(\cS_\ra \rd\big)_\balpha = \sum_{\bbeta\in\ZZ^s}
  \ra_{\balpha-m\bbeta}\; \rd_\bbeta\;,\quad \balpha\in\ZZ^s\;,
  \end{equation}
where $\ell(\ZZ^s)$ is the space of scalar sequences indexed by $\ZZ^s$.
A \emph{subdivision scheme} is the recursive algorithm given by
\begin{equation}\label{eq:subrecursion}
  \rd^{(r+1)}=\cS_\ra \rd^{(r)}=\sum_{\bbeta\in\ZZ^s}
 \ra_{\balpha-m\bbeta}  \rd^{(r)}_\bbeta\;, \quad \rd^{(0)} \in \ell(\ZZ^s), \quad r \in \NN_0\;,
\end{equation}
where $\NN_0$ is the set of natural numbers including zero.

The \emph{symbol} of a subdivision scheme
is given by the Laurent polynomial
\begin{equation*}\label{eq:laurpol}
  a(\bz) = \sum_{\balpha\in\ZZ^s} \ra_\balpha\;\bz^\balpha, \quad \bz = (z_1,\dots,z_s) \in\left( \CC\setminus\{0\}\right)^s,
\end{equation*}
where for $\balpha = (\alpha_1,\dots,\alpha_s) \in \ZZ^s$ we define
$
  \bz^\balpha = z_1^{\alpha_1}z_2^{\alpha_2}\cdot \ldots \cdot z_s^{\alpha_s}.$ Denoting with
\begin{equation}\label{eq:extrpts}
  \rE = \{0,\dots,|m|-1\}^s\,,
\end{equation}
the set of representatives of $\ZZ^s/m  \ZZ^s$ containing $\bnull = (0,0,\dots,0)$,
the $|m|^s$ \emph{submasks} and their symbols $a_\be(\bz)$ are
defined by
\begin{equation} \label{eq:submasks}
   \{\ra_{\be+m\balpha},\ \balpha\in\ZZ^s\}
  \end{equation}
  and
  \begin{equation} \label{eq:subsymbols}
  \quad a_\be(\bz) = \sum_{\balpha\in\ZZ^s}
  \ra_{\be+m\balpha}\;\bz^{\be+m\balpha}
  \;,\quad \be\in\rE\;,
\end{equation}
respectively. Then, we get the  following decomposition of the
mask symbol
\begin{equation} \label{eq:decomp}
  a(\bz) =
  \sum_{\be\in\rE} \; a_{\be}(\bz)\,.
\end{equation}

\subsection{Sum rules and zero conditions}
The \emph{sum rules of order~$1$} in terms of submasks read as follows
\begin{equation}\label{eq:neccond}
  a_\be(\bone) = \sum_{\balpha\in\ZZ^s} \ra_{\be+m\balpha} = 1\;,
  \quad \be\in\rE\;.
\end{equation}
For $\bxi=(\xi_1,\dots,\xi_s) \in \RR^s$ set
$z_j=e^{-i\pi\xi_j}$, $j=1,\dots,s$, the set in \eqref{eq:extrpts} corresponds to

\begin{equation}\label{eq:extrptstwo}
  \Xi=\Xi_\rE = \{\beps =  e^{-i \frac{2 \pi}{m} \be } \;:\; \be \in \rE\}\,,
\end{equation}
and contains $\bone=(1, 1, \dots, 1)$. The sum rules of order~$1$ take an equivalent form
\begin{equation} \label{eq:sumruleone}
 a(\bone) =|m|^s \quad \text{and} \quad a(\beps) = 0 \quad
  \text{for }\beps \in \Xi':=\Xi\setminus \{\bone\}\;.
\end{equation}
Following the notation in \cite{CharinaContiJetterZimm2010}, we call $\Xi'$ the zero set, and the conditions in
\eqref{eq:sumruleone} the \emph{zero condition of order one
(Condition~Z$_1$)}. In the literature, both the conditions in
\eqref{eq:neccond} and their equivalent form
 in \eqref{eq:sumruleone} are   called the sum rules of order one.
We  also make use of the higher order sum rules, see \cite{JePlo} and references therein: The mask symbol
$\ra(\bz)$ is said to satisfy the \emph{zero condition of order
$k$ (Condition~Z$_k$)}, if
\begin{equation} \label{eq:Zkcond}
a(\bone)=|m|^s \quad \hbox{and} \quad
\bigl(D^\bj a\bigr)(\beps)=0 \quad
\hbox{for} \quad \beps\in \Xi'\quad \hbox{and} \quad
|\bj|<k\;.
\end{equation}
We denote by $D^\bj$  the $\bj$-th directional directional
derivative.

\subsection{Parametrization}

\smallskip Since most of the properties of a subdivision scheme, e.g. its convergence, smoothness or its support size, do
not depend on the choice of the associated parameter values
$\bt_\balpha^{(r)}$, $\alpha \in \ZZ^s$, to which the data
$\rd^{(r)}_{\balpha}$, $\alpha \in \ZZ^s$, generated by the $r$-th
step of a subdivision recursion is attached,  these are usually
set to
\begin{equation}\label{def:primal_par}
   \bt_\balpha^{(r)}:=\frac{\balpha}{m^r}, \quad \balpha \in \ZZ^s, \quad r\ge 0.
\end{equation}
We refer to the choice in (\ref{def:primal_par}) as
\emph{standard} parametrization. We show in section
\ref{sec:algebra} that the capability of subdivision to reproduce
polynomials  does depend on the choice of the associated parameter
values and the standard parametrization is not always the optimal
one. As in \cite{ContiHormann} for the univariate case, the choice
\begin{equation}\label{def:general_par}
   \bt_\balpha^{(r)}:= \bt_\bnull^{(r)}+\frac{\balpha}{m^r},
   \quad \bt_\bnull^{(r)}=\bt_\bnull^{(r-1)}-\frac{\btau}{m^r}, \quad \bt_\bnull^{(0)}=0,\quad \balpha \in \ZZ^s, \quad r\ge 0,
\end{equation}
with a suitable $\btau \in \RR^s$ turns out to be a better selection.\\
We call the sequence $\{\bt^{(r)},\ k\ge 0\}$, with
$\bt^{(r)}=\{\bt^{(r)}_\balpha,\ \alpha\in \ZZ^s\}$  the \emph{sequence
of parameter values associated with the subdivision scheme}.

\subsection{Convergence and non-singularity of subdivision}
Following \cite{DynShen2008}, our definition of convergence
depends on the parameter values associated with a given
subdivision scheme. Since the subdivision process generates denser and denser sequences of data $\rd^{(r)}$, $r \ge 0$, a notion of
convergence can be established by using a sequence $\{F^{(r)}, r \ge 0 \}$ of continuous functions $F^{(r)}$
that interpolate the data $\rd^{(r)}$ at the parameter values $\{\bt^{(r)},\ k\ge 0\}$ associated to the subdivision scheme,
namely
\begin{equation} \label{eq:interpolating_functions}
  F^{(r)}( \bt_\balpha^{(r)}) = \rd_\balpha^{(r)}, \qquad \balpha \in \ZZ^s,\quad  r\geq0.
\end{equation}

\begin{definition}
If the sequence of continuous functions $\{F^{(r)},\ r\ge 0\}$ satisfying \eqref{eq:interpolating_functions}
converges, then we denote its limit by
\[
  g_\rd := \lim_{r \to\infty} F^{(r)},
\]
and say that $g_\rd$ is the \emph{limit function} of the associated subdivision scheme (\ref{eq:subrecursion}) for
the initial data $\rd^{(0)}=\{\rd^{(0)}_\balpha,\ \balpha\in \ZZ^s\} \in \ell(\ZZ^s)$. The limit function $\phi:=g_\bdelta$
with the initial data
$$
  \delta_\alpha=\left\{ \begin{array}{cc} 1, & \alpha=\bnull, \\ 0, & \hbox{otherwise}\end{array}\right. , \quad
\alpha \in \ZZ^s,
$$
is called the basic limit function of this scheme, it is compactly supported and satisfies the refinement equation

\begin{equation} \label{refinement_equation}
 \phi=\sum_{\balpha \in \ZZ^s} \ra_\balpha \phi(m\cdot-\balpha)
\end{equation}
with refinement coefficients given by the mask $\ra$.
\end{definition}


\begin{definition}
A subdivision scheme is called \emph{non-singular}, if it is convergent, and $g_\rd=0$ if and only if $\rd$ is the zero sequence, i.e. $\rd_\balpha=0$ for
all $\balpha\in \ZZ^s$.
\end{definition}

\smallskip 
Next we show that
the notion of non-singular subdivision scheme is equivalent to the notion of linear independence of the integer shifts of its basic limit function.

\begin{proposition} \label{prop:nonsingular_vs_independence}
A convergent subdivision scheme $S_\ra$ is non-singular if and
only if the integer translates of the solution $\phi$ of
refinement equation (\ref{refinement_equation}) with coefficients satisfying
$a_\be(\bone)=1,\ \be\in \rE$, are linearly independent and form a
partition of unity.
\end{proposition}
\pf

\noindent ''$\Longrightarrow$:`` Assume that $S_\ra$ is convergent, then
the basic limit function  $\phi= S_\ra^{\infty} \bdelta$ satisfies the refinement equation \eqref{refinement_equation} and its
integer shifts form a partition of unity.
If the convergent subdivision scheme $S_\ra$ is non-singular, then
for any starting sequence $\rd \in \ell(\ZZ^s)$ we have
$$
 S_\ra^{\infty} \rd=\sum_{\balpha \in \ZZ^s} \rd_\balpha \phi(\cdot-\balpha)=0
$$
if and only if  $\rd$ is the zero sequence.

\noindent ''$\Longleftarrow$:`` Assume that the integer translates of
the solution $\phi$ of refinement equation \eqref{refinement_equation}
are linearly independent and form a partition of unity. This implies that the subdivision scheme $S_\ra$ associated with the
symbol $a(\bz)$ is convergent, see \cite[Lemma 2.3]{D92}. It is also then non-singular, otherwise one easily gets a
contradiction to the assumption on linear independence of the translates of $\phi$.
\eop

\subsection{Polynomial generation versus polynomial reproduction}
We denote by $\Pi_k$ the space of multivariate  polynomial of total degree $k \in \NN_0$.

\begin{definition}[Polynomial generation] \label{def:pg}
A convergent stationary subdivision scheme $S_\ra$ \emph{generates} polynomials up to degree
$d_G$ (is $\Pi_{d_G}$-generating) if for any polynomial $\pi \in \Pi_{d_G}$ there exists some
initial data ${\rrq}^{(0)} \in \ell_\infty(\ZZ^s)$ such that $S^{\infty}_{\ra}{\rrq}^{(0)}=\pi$.
Moreover, the initial data ${\rrq}^{(0)}$ is sampled from a polynomial $\widetilde{\pi} \in \Pi_{d_G}$
 with the same leading coefficients as $\pi \in \Pi_{d_G}$.
\end{definition}
Note that the assumptions on the properties of $\widetilde{\pi}$ in the above Definition are justified
by \cite[Lemma 2.1]{ALevin2003}. Note also that polynomial generation is also studied in \cite{CDM91}.\\

\noindent We continue with a slightly different notion, the notion of polynomial reproduction which requires a specific choice of starting sequences of a polynomial limit. The concepts of polynomial reproduction and
generation coincide for $d_G=d_R=0$.

\begin{definition}[Polynomial reproduction] \label{def:pr}
A convergent subdivision scheme $S_{{\ra}}$ with parameter values $\{\bt^{(r)},\ r\ge 0\}$
is \emph{reproducing} polynomials up to degree $d_R$
(is $\Pi_{d_R}$-reproducing) if for any polynomial $\pi \in \Pi_{d_R}$ and for the initial data
${\rp}^{(0)}=\{\pi(\bt^{(0)}_\balpha),\ \balpha\in\ZZ^s\}$ the limit of the subdivision satisfies
$S^{\infty}_{\ra}{\rp}^{(0)}=\pi$.
\end{definition}

Another important property of subdivision is the so-called step-wise polynomial reproduction, we make use of it in Section 2.

\begin{definition} [Step-wise polynomial reproduction] \label{def:swpr}
A convergent subdivision scheme $S_\ra$ with parameter values $\{\bt^{(r)},\ r\ge 0\}$ is  \emph{step-wise polynomial reproducing}  up
to degree $k$ (is step-wise $\Pi_{k}$-reproducing) if for any polynomial $\pi \in \Pi_k$ and for the
data $\rd^{(r)}=\{\pi(\bt^{(r)}_\balpha),\ \balpha\in\ZZ^s\}$
\begin{equation}\label{passo}
    \rd^{(r+1)}=S_{\ra}\rd^{(r)}\quad \hbox{or, equivalently,} \quad
 \pi(\bt^{(r+1)}_\balpha)=\sum_{\bbeta\in \ZZ^s}\ra_{\balpha-m\bbeta} \, \pi(\bt^{(r)}_\bbeta),\ \ \balpha\in \ZZ^s\,.
\end{equation}
\end{definition}

The next proposition shows that for a non-singular subdivision scheme the concepts of polynomial reproduction and
step-wise polynomial reproduction are equivalent.

\begin{proposition}\label{step-limit}
A non-singular, subdivision scheme $S_\ra$  is \emph{step-wise polynomial reproducing} up
to degree $k$ if and only if it is \emph{polynomial reproducing} up  to degree $k$.
\end{proposition}

\pf \vspace{0.1cm} \noindent ''$\Longrightarrow$:`` For any polynomial $\pi \in \Pi_{k}$, let ${\rd}^{(0)}:=\{\pi(\bt_\balpha^{(0)}),\ \balpha\in \ZZ^s\}$.
If the subdivision scheme $S_\ra$ is step-wise $\Pi_{k}$-reproducing, then the sequence
$\{F^{(r)}, r \ge 0\}$ of continuous functions $F^{(r)}$ satisfying \eqref{eq:interpolating_functions}
with  $\rd_\balpha^{(r)}=\pi\left(\bt_\balpha^{(r)}\right)$, $\alpha \in \ZZ^s$,
converges uniformly to $\pi$ as $r \rightarrow\infty$, since the distance between the grid  points $m^{-r}\ZZ^s$ goes to zero.

\vspace{0.1cm} \noindent ''$\Longleftarrow$:`` Let us assume next that the subdivision scheme $S_\ra$ is \emph{$\Pi_{k}$-reproducing}. Let $r \ge 0$.
On the one hand, applying  the subdivision scheme to the data ${\rd}^{(r)}=\{\pi(\bt_\balpha^{(r)}),\ \balpha\in \ZZ^s\}$ we obtain,
$$
S^\infty_\ra{\rd}^{(r)}=S^\infty_\ra \rd^{(r+1)}=\pi, \quad \rd^{(r+1)}=S_\ra \rd^{(r)}.
$$
On the other hand, for the sequence ${\rp}^{(r+1)}:=\{\pi(\bt_\balpha^{(r+1)}),\ \balpha\in \ZZ^s\}$ we also have
$$
S^\infty_{\ra}{\rp}^{(r+1)}=\pi\,.
$$
Therefore, by the linearity of the operator $S_{a}$ it follows
$$
S^\infty_{\ra}\left({\rd}^{(r+1)}-{\rp}^{(r+1)}\right)=0
$$
and, thus, ${\rp}^{(r+1)}=\rd^{(r+1)}=S_{\ra}{\rd}^{(r)}$ due to the assumption of non-singularity. Thus,
the claim follows.
\eop

\section{Algebraic condition for polynomial generation and reproduction} \label{sec:algebra}
In this section, for a non-singular subdivision scheme, we
determine the value of $\btau \in \RR^s$ in
\eqref{def:general_par} that guarantees the polynomial
reproduction of linear polynomials, see Proposition
\ref{prop:linear-reproduction}.  In Theorem \ref{teo:main}, we
then provide algebraic conditions on $a(\bz)$  for checking the
reproduction of polynomials of higher degree. These algebraic
conditions depend on the previously obtained value of $\btau$. In
the case of only convergent schemes see Corollaries
\ref{cor:linear-reproduction} and \ref{cor:main}. We start by
defining  the tensor product polynomial of degree $|\bj|$, $\bj \in \NN_0^s$, given by
 \begin{equation}\label{def:qj}
q_{\bnull}(\bz):=1,\quad  q_{\bj}(z_1,\dots, z_s):=\prod_{i=1}^s
\prod_{\ell_i=0}^{j_i-1}(z_i-\ell_i), \quad \bj=(j_1,\dots, j_s).
 \end{equation}
The following auxiliary proposition states known results on polynomial generation
that we make use of in this section.

\begin{proposition}\label{prop:equiv}
Let $S_\ra$ be a convergent subdivision scheme.

\begin{description}
 \item[$(i)$] The subdivision scheme $S_\ra$ reproduces constant sequences or, equivalently, its
symbol $a(\bz)$ satisfies
  (\ref{eq:neccond}), if and only if $a(\bz)$ satisfies \eqref{eq:sumruleone}.
 \item[$(ii)$]  Let $k \in \NN$. The  symbol $a(\bz)$ satisfies condition  $Z_k$ if and only if
  $$\bigl(D^\bj a_\be\bigr)(\bone)=|m|^{-s} D^\bj a(\bone),\quad
  \hbox{for} \quad  \be\in \rE, \quad \bj \in \NN_0^s, \quad |\bj|< k.$$
 \item[$(iii)$]  Let $k \in \NN$.  The symbol $a(\bz)$ satisfies Condition  $Z_k$  if and only if
 \begin{equation}\label{eq:equiv2}
 \sum_{\bbeta\in \ZZ^2}q_{\bj}(\balpha-m\bbeta) \ra_{\balpha-m\bbeta}=|m|^{-s} D^{\bj}a(\bone),\quad \alpha\in \ZZ^s,
\quad \bj \in \NN_0^s, \quad |\bj|<k.
 \end{equation}
\end{description}
\end{proposition}

\pf The proof of $(i)$ for dilation matrix $2I$ follows from the
definition of the operator $S_\ra$ and \cite[Section 6]{CDM91} and
in the case of dilation matrix $mI$, $|m|>2$, is in \cite[Lemma 3.3
]{Jia98}. The proof of $(ii)$: Let $\bj \in \NN_0^s$, $|\bj|< k$. Since
$a_\be(\bz)=\displaystyle{\sum_{\balpha\in\ZZ^s}
\ra_{\be-m\balpha}\bz^{\be-m\balpha}}$, its  $\bj$-th derivative
satisfies
\begin{equation}\label{eq:dersubsymbol}
\bigl(D^\bj
a_\be\bigr)(\bz)=\sum_{\balpha\in\ZZ^s}q_{\bj}(\be-m\balpha)
\ra_{\be-m\balpha}\bz^{\be-m\balpha-\bj}, \quad \bz \in (\CC
\setminus \{ 0\})^s.
\end{equation}
Next, due to $\displaystyle a(\bz)=\sum_{\be \in \rE}
a_{\be}(\bz)$, we have
$$
\begin{array}{ll}
  \bigl(D^\bj a\bigr)(\beps)&=\displaystyle{\sum_{\be \in \rE} \bigl(D^\bj a_\be\bigr)(\beps)}=
   \displaystyle{\sum_{\be \in \rE}\, \sum_{\balpha\in\ZZ^s}q_\bj (\be-m\balpha) \ra_{\be-m\balpha}\beps^{\be-m\balpha-\bj}}\\
  \\
  & =\displaystyle{\sum_{\be \in \rE}\beps^{\be-\bj}\, \sum_{\balpha\in\ZZ^s}q_\bj (\be-m\balpha) \ra_{\be-m\balpha}}=
 \displaystyle{\sum_{\be \in \rE}\beps^{\be-\bj}\, \bigl(D^\bj a_\be\bigr)(\bone)}\\
\end{array}
$$
for $\beps \in \Xi$. The rows of the  matrix of this linear system
are given by $(\beps^{\be-\bj})_{\be \in \rE}$. This matrix is
invertible and Condition $Z_k$, due to
$$
 \sum_{\be \in \rE} \beps^{\be-\bj}=\left\{ \begin{array}{cc} |m|^s, & \beps=1, \\ 0, & \hbox{otherwise}, \end{array} \right.
$$
is equivalent to
$$
 \sum_{\be \in \rE} D^\bj a_\be(\bone)= D^\bj a(\bone), \quad  D^\bj a_\be(\bone)= D^\bj a_{\widetilde{\be}}(\bone),
\quad \be, \widetilde{\be} \in \rE, \quad \be \not =
\widetilde{\be}.
$$
Proof of $(iii)$:  
For any $\balpha\in \ZZ^s$ there exists $\be\in \rE$ and $\bbeta
\in \ZZ^s$ such that $\balpha=\be+m\beta$. Thus,
$$
 D^{\bj}a_\be(\bone)=\sum_{\bbeta\in\ZZ^s}q_{\bj}(\balpha-m\bbeta)\ra_{\balpha -m\bbeta},
\quad \be \in E\,.
$$
The claim follows by $(ii)$.
\eop

\begin{remark} Note that the proof of $(ii)$ is also implied, for example, by \cite[Theorem 3.7]{Cabrelli}.
\end{remark}

The next result provides a  simple algebraic conditions for determining $\btau \in \RR^s$, which appears in \eqref{def:general_par}
and guarantees the reproduction of linear polynomials.

\begin{proposition}\label{prop:linear-reproduction}
Let $S_\ra$ be a non-singular subdivision scheme that generates linear
polynomials, i.e. its symbol satisfies Condition $Z_1$. Then $S_\ra$ reproduces linear
polynomials if and only if its parameter values are given by (\ref{def:general_par})
with $\btau=|m|^{-s} \left(D^{\bveps_1} a(\bone),\dots, D^{\bveps_s} a(\bone)\right)$.
\end{proposition}

\pf According to Proposition \ref{step-limit} for non-singular subdivision schemes, polynomial reproduction
is equivalent to step-wise polynomial reproduction. Moreover, any convergent subdivision scheme reproduces
the constants, hence it is sufficient to prove the claim for polynomials of the form $\pi(x_1, \dots, x_s)=x_j$,
$j=1, \dots,s$. Let $r\in\NN_0$ and set  $\rd^{(r)}_\balpha=\pi(\bt^{(r)}_\balpha)$, $\balpha \in \ZZ^s$, with $\pi(\bx)=x_j$.
Then  for any $\balpha\in\ZZ^s$ and $\be=(\be_1, \dots, \be_s)$ we get
\begin{align*}
  \rd^{(r+1)}_{m\balpha+\be}
  &= \sum_{\bbeta\in\ZZ^s} \ra_{m(\balpha-\bbeta)+\be} \rd^{(r)}_\bbeta
   = \sum_{\bbeta\in\ZZ^s} \ra_{m\bbeta+\be} \rd^{(r)}_{\balpha-\bbeta}
   = \sum_{\bbeta\in\ZZ^s} \ra_{m\bbeta+\be} \left(t_{\bnull,j}^{(r)}+\frac{\alpha_j-\beta_j}{m^r}\right)\\
  &= \underbrace{\sum_{\bbeta\in\ZZ^s} \ra_{m\bbeta+\be}}_{a_\be(\bone)} \left(t_{\bnull,j}^{(r)}+\frac{m \alpha_j+e_j}{m^{r+1}}\right)
     - \sum_{\bbeta\in\ZZ^s} \ra_{m\bbeta+\be} \frac{m \beta_j+e_j}{m^{r+1}}\\
  &= \left(t_{\bnull,j}^{(r)}+\frac{m \alpha_j+e_j}{m^{r+1}}\right) - \frac{D^{\bveps_j} a_\be(\bone)}{m^{r+1}}\\
   &= \left(t_{\bnull,j}^{(r)}+\frac{m \alpha_j+e_j}{m^{r+1}}\right) - \frac{D^{\bveps_j} a(\bone)}{ |m|^s \cdot m^{r+1}},
\end{align*}
where the last equality is due to Proposition \ref{prop:equiv} part $(ii)$ for $\bj=\bveps_j$. Thus, $ \rd^{(r+1)}_{m\balpha+\be}$ is equal to
\[
  \pi(t^{(r+1)}_{m\balpha+\be})=t_{\bnull,j}^{(r+1)}+\frac{m \alpha_j+e_j}{m^{r+1}}=
  t_{\bnull,j}^{(r)}+\frac{m \alpha_j-\tau_j+e_j}{m^{r+1}}, \quad \alpha \in \ZZ^s,
\]
if and only if $\tau_j=|m|^{-s} D^{\bveps_j} a(\bone)$.
\eop

Note that not all convergent subdivision schemes are non-singular, e.g. the subdivision  scheme
based on the four directional box spline symbol is not. For such schemes one can still
determine a parametrization that ensures its polynomial generation property. The following
result is a direct consequence of Proposition \ref{prop:linear-reproduction}.

\begin{corollary} \label{cor:linear-reproduction}
Let $S_\ra$ be a convergent subdivision scheme that generates linear
polynomials, i.e. its symbol satisfies Condition $Z_1$. Then $S_\ra$ reproduces linear
polynomials if its parameter values are given by (\ref{def:general_par})
with $$\btau=|m|^{-s} \left(D^{\bveps_1} a(\bone),\dots, D^{\bveps_s} a(\bone)\right).$$
\end{corollary}

The following result is crucial for the proof of the main result, Theorem \ref{teo:main} of this section. It allows us to express the
polynomial generation of $S_\ra$ in terms of the properties of its symbol.

\begin{proposition}\label{prop:equiv3}
Let $k \in \NN$, $\btau\in \RR^s$ and $q_{\bj}$ as in
\eqref{def:qj}. A subdivision symbol $a(\bz)$ satisfies
\begin{equation}\label{eq:equiv1a}
 \bigl(D^\bj a\bigr)(\bone)=|m|^{s} q_{\bj}(\btau), \quad \bigl(D^\bj a\bigr)(\beps)=0
 \quad    \beps\in \Xi',\quad \bj \in \NN_0^s, \quad |\bj| \le k,
\end{equation}
if and only if
\begin{equation}\label{eq:equiv2a}
\sum_{\bbeta\in \ZZ^s} \ra_{\balpha-m\bbeta} \bbeta^\bj
=\left(\frac{\balpha-\btau}{m}\right)^\bj, \quad \alpha\in \ZZ^s,
\quad \bj \in \NN_0^s, \quad |\bj|\le k\,.
\end{equation}
\end{proposition}
\pf  First, note that due to Proposition \ref{prop:equiv} the conditions in \eqref{eq:equiv1a}
are equivalent to
\begin{equation}\label{eq:propqj}
  q_{\bj}(\btau)=\sum_{\bbeta\in \ZZ^s} q_{\bj}(\balpha-m\bbeta)\ra_{\balpha-m\bbeta}, \quad
\bj \in \NN_0^s, \quad |\bj| \le k\,.
\end{equation}
Using this equivalent formulation, we prove the proposition by induction on $k$. For $k=0$ the claim is true since
for any $\btau \in \RR^s$ we get
$$
  q_{\bnull}(\btau)=\sum_{\bbeta\in \ZZ^s}\ra_{\balpha-m\bbeta}=\left(\frac{\balpha-\btau}{m}\right)^{\bnull}=1.
$$
Next, we assume that the claim is true for all $|\bj| \le k-1$ and prove it for $\bj \in \NN_0^s$ with $|\bj|=k$.
To this purpose,  we write the polynomial $q_{\bj}$ in $\bx$ of (total) degree $|\bj|=k$ as
\begin{equation}\label{eq:propqj2}
q_{\bj}(\balpha-m\bx)=\sum_{\bell \in \NN_0^s, \ |\bell|\le k} c_{\bj,\balpha,\bell} \bx^\bell,
\quad \bx \in \RR^s,\quad  c_{\bj,\balpha,\bj}\neq 0\,.
\end{equation}
 Therefore, using the induction assumption  and by \eqref{eq:propqj} and \eqref{eq:propqj2} we have
$$\begin{array}{ll}
    q_{\bj}(\btau)&=\displaystyle{\sum_{\bbeta\in \ZZ^s}\sum_{\bell  \in \NN_0^s, \ |\bell|\le k}
    c_{\bj,\balpha,\bell}  \ra_{\balpha-m\bbeta}} \bbeta^\bell \\\\
    &=\displaystyle{\sum_{\bell  \in \NN_0^s, \ |\bell|= k}c_{\bj,\balpha,\bell} \sum_{\bbeta\in \ZZ^s}
    \ra_{\balpha-m\bbeta} \bbeta^\bell}
    +\displaystyle{\sum_{\bell  \in \NN_0^s, \ |\bell|\le k-1}c_{\bj,\balpha,\bell} \sum_{\bbeta\in \ZZ^s}
     \ra_{\balpha-m\bbeta} \bbeta^\bell}\\
    \\
     &=\displaystyle{\sum_{\bell  \in \NN_0^s, \ |\bell|= k}c_{\bj,\balpha,\bell} \sum_{\bbeta\in \ZZ^s}
     \ra_{\balpha-m\bbeta} \bbeta^\bell}+\displaystyle{\sum_{\bell  \in \NN_0^s, \ |\bell|\le k-1}c_{\bj,\balpha,\bell} \left(\frac{\balpha-\btau}{m}\right)^\bell}\\
     \\
      &=\displaystyle{\sum_{\bell  \in \NN_0^s, \ |\bell|= k}c_{\bj,\balpha,\bell} \left(\sum_{\bbeta\in \ZZ^s}
       \ra_{\balpha-m\bbeta} \bbeta^\bell-
      \left(\frac{\balpha-\btau}{m}\right)^\bell\right)}+q_{\bj}(\btau)\,.
  \end{array}
    $$
    The last equality is due to the fact that
    $$
      q_\bj(\btau)=q_\bj\left(\balpha-m \cdot \frac{\balpha-\btau}{m} \right)=\sum_{\bell  \in \NN_0^s, \ |\bell|\le k} c_{\bj,\balpha,\bell} \left(\frac{\balpha-\btau}{m}\right)^\bell.
    $$
    Hence, due to  $c_{\bj,\balpha, \bj} \not =0$, the above identity holds if and only if
$$
\sum_{\bbeta\in \ZZ^s} \ra_{\balpha-m\bbeta} \bbeta^\bj-
      \left(\frac{\balpha-\btau}{m}\right)^\bj=0,\quad \hbox{for}\  \bj \in \NN_0^s, \quad
      |\bj|=k,
$$
and the claim follows.\eop

\smallskip \noindent
We are now ready to prove the main results of this paper.
\begin{theorem}\label{teo:main}
Let $k\in \NN_0$. A non-singular subdivision scheme with symbol $a(\bz)$ and associated parametrization in (\ref{def:general_par}) with
some $\btau \in \RR^s$ reproduces polynomials of degree up to $k$
if and only if
$$
 \bigl(D^\bj a\bigr)(\bone)= |m|^s q_{\bj}(\btau) \quad \hbox{and} \quad \bigl(D^\bj a\bigr)(\beps)=0 \quad
\hbox{for} \quad \beps\in \Xi',\quad |\bj| \le k\,.$$
\end{theorem}
\pf The proof is by induction on $k$. In the case $k=0$ the claim follows by part $(i)$ of Proposition \ref{prop:equiv}. By Proposition \ref{step-limit} it suffices to prove the result for the stepwise polynomial reproduction.

\vspace{0.1cm} \noindent ''$\Longleftarrow$:`` We write any polynomial $\pi$ of degree $k$ as
$\pi(\bx)=\displaystyle{\sum_{\bell \in \NN_0^s, |\bell|= k} c_\bell \bx^\bell + \widetilde{\pi}(\bx)}$ with
$\widetilde{\pi} \in \Pi_{k-1}$.  Let $r \ge 0$. We show that the sequence
$$
 \rd^{(r)}=\left\{\pi(\bt_\balpha^{(r)})= \displaystyle{\sum_{\bell \in \NN_0^s, |\bell|= k}
 c_\bell \left(\bt_\balpha^{(r)}\right)^\bell+ \widetilde{\pi}(\bt_\balpha^{(r)})},\ \balpha\in \ZZ^s\right \}
$$
satisfies  $\rd^{(r+1)}=S_\ra \rd^{(r)}=\{\pi(\bt_\balpha^{(r+1)}),\ \balpha\in \ZZ^s\}$. In fact, due to the
induction assumption, by \eqref{def:general_par} and  Proposition \ref{prop:equiv3}, we have
$$
\begin{array}{ll}
\rd^{(r+1)}_\balpha&=\displaystyle{\sum_{\bbeta\in \ZZ^s}\ra_{\balpha-m\bbeta}\rd^{(r)}_\bbeta}=
\displaystyle{\sum_{\bbeta\in \ZZ^s}\ra_{\balpha-m\bbeta}\sum_{\bell \in \NN_0^s, |\bell|= k} c_\bell
\left(\bt_\bnull^{(r)}+\frac{\bbeta}{m^r}\right)^\bell+\widetilde{\pi}(\bt^{(r+1)}_\balpha)}\\
\\
&=\displaystyle{\sum_{\bbeta\in \ZZ^s}\ra_{\balpha-m\bbeta}\sum_{\bell \in \NN_0^s, |\bell|= k}c_\bell \sum_{\bh\le \bell}\left(
                             \begin{array}{c}
                             \bell \\
                             \bh \\
                        \end{array}
                        \right)
\left(\frac{\bbeta}{m^r}\right)^\bh \left(\bt_\bnull^{(r)}\right)^{\bh-\bell}+\widetilde{\pi}(\bt^{(r+1)}_\balpha)}\\
\\
&=\displaystyle{\sum_{\bell \in \NN_0^s, |\bell|= k}c_\bell \sum_{\bh \le \bell}\left(
                                                                                                                  \begin{array}{c}
                                                                                                                    \bell \\
                                                                                                                    \bh \\
                                                                                                                  \end{array}
                                                                                                                \right)
\left(\bt_\bnull^{(r)}\right)^{\bh-\bell}\left(\frac{1}{m^r}\right)^\bh\sum_{\bbeta\in \ZZ^s}\bbeta^\bh
\ra_{\balpha-m\bbeta}+\widetilde{\pi}(\bt^{(r+1)}_\balpha)}\\ \\
&=\displaystyle{\sum_{\bell \in \NN_0^s, |\bell|= k}c_\bell \sum_{\bh\le \bell}\left(
                                                                                                                  \begin{array}{c}
                                                                                                                    \bell \\
                                                                                                                    \bh \\
                                                                                                                  \end{array}
                                                                                                                \right)
\left(\bt_\bnull^{(r)}\right)^{\bh-\bell}\left(\frac{\balpha-\btau}{m^{r+1}}\right)^\bh+\widetilde{\pi}(\bt^{(r+1)}_\balpha)}\\
\\
&=\displaystyle{\sum_{\bell \in \NN_0^s, |\bell|= k}c_\bell \left(\bt_\balpha^{(r+1)}\right)^\bell}+\widetilde{\pi}(\bt^{(r+1)}_\balpha)=\pi(\bt^{(r+1)}_\balpha),
\quad \balpha \in \ZZ^s.
\end{array}
$$
The one but last equality is due to
$$
\bt_\bnull^{(r)}+\frac{\balpha-\btau}{m^{r+1}}=\bt_\bnull^{(r+1)}+\frac{\balpha}{m^{r+1}}= \bt_\balpha^{(r+1)}\,.
$$

\vspace{0.1cm} \noindent ''$\Longrightarrow$:`` Let $\bj \in
\NN_0^s$ be such that $|\bj|=k$ with $j_i=k$ for some $i=1, \dots,
s$. Let the polynomial $\pi(\bx)=\bx^\bj$ and the sequence
$\rd^{(r)}=\{\pi(\bt_\balpha^{(r)}),\ \balpha\in \ZZ^s\}$. On the
one hand,  by similar arguments as above, we get
$$
  \rd^{(r+1)}_\balpha= \sum_{\bh  \le \bj }\left(
                                                                                                                  \begin{array}{c}
                                                                                                                    \bj \\
                                                                                                                    \bh \\
                                                                                                                  \end{array}
                                                                                                                \right)
\left(\bt_\bnull^{(r)}\right)^{\bh-\bj}\left(\frac{1}{m^r}\right)^\bh \sum_{\bbeta\in \ZZ^s}\bbeta^\bh
\ra_{\balpha-m\bbeta}.
$$
On the other hand, the definition of $\bt_\bnull^{(r+1)}$ yields
$$
 \pi(\bt^{(r+1)}_\balpha)= \sum_{ \bh \le \bj}\left(
                                                                                                                  \begin{array}{c}
                                                                                                                    \bj \\
                                                                                                                    \bh \\
                                                                                                                  \end{array}
                                                                                                                \right)
\left(\bt_\bnull^{(r)}\right)^{\bh-\bj} \left(\frac{\balpha-\btau}{m^{r+1}}\right)^\bh.
$$
The polynomial reproduction, \ie $\ \rd^{(r+1)}=S_\ra \rd^{(r)}$,
implies that
 $$
   \displaystyle{\sum_{\bh \le \bj}\left(
                                                                                                                  \begin{array}{c}
                                                                                                                    \bj \\
                                                                                                                    \bh \\
                                                                                                                  \end{array}
                                                                                                                \right)
\left(\bt_\bnull^{(r)}\right)^{\bh-\bj}\left(\frac{1}{m^r}\right)^\bh \left[\sum_{\bbeta\in \ZZ^s}\bbeta^\bh
\ra_{\balpha-m\bbeta}-\left(\frac{\balpha-\btau}{m}\right)^\bh\right]}=0.
$$
Thus, by induction for $\bh=\bj$ we have
$$
  \displaystyle \sum_{\bbeta\in \ZZ^s} \bbeta^\bj
\ra_{\balpha-m\bbeta}-\left(\frac{\balpha-\btau}{m }\right)^\bj=0, \quad \alpha \in \ZZ^s.
$$
The claim follows from Proposition \ref{prop:equiv3}.
\eop

For convergent schemes we readily get the following result, which is due to the fact that for
convergent schemes the
step-wise polynomial reproduction implies polynomial reproduction.

\begin{corollary}\label{cor:main}
Let $k\in \NN_0$. A convergent subdivision scheme with symbol $a(\bz)$ and associated parametrization in (\ref{def:general_par}) with
some $\btau \in \RR^s$ reproduces polynomials of degree up to $k$
if
$$
 \bigl(D^\bj a\bigr)(\bone)= |m|^s q_{\bj}(\btau) \quad \hbox{and} \quad \bigl(D^\bj a\bigr)(\beps)=0 \quad
\hbox{for} \quad \beps\in \Xi',\quad |\bj| \le k\,.$$
\end{corollary}

\section{Applications and examples}

It is natural to expect that  any shift of the subdivision mask does not effect the polynomial reproduction properties
of the corresponding scheme, which is confirmed by the next result.

\begin{lemma} A convergent subdivision scheme $S_\ra$ with the symbol $a(\bz)$ reproduces polynomials up to degree $k$ if and only if
so does the shifted scheme $S_{\widetilde{\ra}}$ with the symbol $\widetilde{a}(\bz)=\bz^\balpha a(\bz)$, $\balpha \in \ZZ^s$.
\end{lemma}
\pf By Proposition \ref{prop:linear-reproduction} and due to
$a(\bone)=m^s$, we get the following identity for the suitable
$\btau$ of $S_\ra$ and $\widetilde{\btau}$ of
$S_{\widetilde{\ra}}$
$$
  \widetilde{\btau}=\btau+\balpha.
$$
By Leibnitz differentiation formula and due to the fact that $D^\bell \bz^\balpha$, $\bell \in \NN_0^s$, evaluated at
$\bone$ is equal to  $q_\bell(\balpha)$ in \eqref{def:qj} we have
$$
 D^{\bj} \widetilde{a}(\bone)=\sum_{\bell \in \NN_0^s, \bell \le \bj} \left( \begin{array}{c} \bj \\ \bell \\ \end{array} \right)
 q_\bell(\balpha)  D^{\bj-\bell} a(\bone).
$$
Thus, by Corollary \ref{cor:main}, to prove the claim we need to
show that
$$
  q_\bj(\btau+\balpha)=\sum_{\bell \in \NN_0^s, \bell \le \bj} \left( \begin{array}{c} \bj \\ \bell \\ \end{array} \right)
 q_\bell(\balpha)  q_{\bj-\bell}(\btau).
$$
By definition of $q_\bj$ it is a tensor product polynomial, thus, it suffices to show that the following two
univariate polynomials in $\tau \in \RR$ are equal
$$
 \prod_{n=0}^{j-1}(\tau+\alpha-n) =\sum_{0 \le \ell \le j} \left( \begin{array}{c} j \\ \ell \\ \end{array} \right)
  \prod_{i=0}^{\ell-1}(\tau-i)  \prod_{t=0}^{j-\ell-1}(\alpha-t).
$$
The claim follows by the one dimensional result \cite[Corollary 5.1]{ContiHormann}.
\eop

\subsection{Box splines}
An $s-$variate box spline is given by its symbol
\begin{equation} \label{def:box_spline_symbol}
  a_\Theta(\bz)=2^s \prod_{\theta \in \Theta} \frac{1+\bz^\theta}{2},
\end{equation}
where $\theta$ runs through all the columns of the $s \times n$,
 rank $s$ matrix $\Theta\in \ZZ^{s\times n}$ with $n \ge s$.
It is well-known that the subdivision schemes associated with the
symbols $a_\Theta(\bz)$ are convergent, if the matrix $\Theta$
is such that removing any column from $\Theta$ does not change its
rank, see \cite[p. 127]{BHR}. Next results gives the correct parametrization  for box spline subdivision schemes.

\begin{lemma}\label{lemma:box}
A subdivision scheme with the symbol $\ra_\Theta(\bz)$ in \eqref{def:box_spline_symbol} reproduces
linear polynomials if its associated parameter values are as in \eqref{def:general_par}  with
$$
 \btau=\frac{1}{2} \left(\sum_{\theta \in \Theta} \theta_1, \dots, \sum_{\theta \in \Theta} \theta_s \right),
 \quad \theta=(\theta_1, \dots, \theta_s).
$$
\end{lemma}
\pf
 The result follows from Corollary \ref{cor:linear-reproduction} and the simple fact that
  $$
   D^{\epsilon_j} a_\Theta(\bz)=2^s \cdot \frac{1}{2} \cdot \sum_{\theta \in \Theta} \theta_j
   z^{\theta-\epsilon_j}
   \prod_{\widetilde{\theta} \in \Theta \atop \widetilde{\theta} \not =\theta} \frac{1+\bz^{\widetilde{\theta}}}{2},\quad z\in (\CC\setminus \{0\})^s.
  $$
\eop

\begin{remark} In case $\Theta$ is unimodular, i.e. each $s \times s$ submatrix of $\Theta$ has determinant  $1$ or
$-1$, the integer shifts of the corresponding box splines are
linear independent and therefore the subdivision scheme associated
with $\ra_\Theta(\bz)$, if convergent,  is non-singular due to
Proposition \ref{prop:nonsingular_vs_independence}. Hence, the
results of Proposition \ref{prop:linear-reproduction} and of
Theorem  \ref{teo:main} hold.
\end{remark}

\medskip We consider an example of the $3-$directional box splines with
$$
 \Theta=\left( \begin{array}{ccc} \underbrace{\begin{array}{ccc} 1& \dots & 1
 \\ 0& \dots & 0 \end{array}}_{\hbox{$k$ times}} & \underbrace{\begin{array}{ccc} 0 & \dots & 0 \\
   1 & \dots & 1 \end{array}}_{\hbox{$\ell$ times}} & \underbrace{\begin{array}{ccc}  1 & \dots & 1 \\
   1& \dots & 1 \end{array}}_{\hbox{$n$ times}} \end{array} \right)
$$
and the corresponding symbols $a_\Theta(z)$ are denoted by
$$
 B_{k, \ell, n}(z_1,z_2)= 4 \cdot  \left( \frac{1+z_1}{2} \right)^k \left( \frac{1+z_2}{2} \right)^\ell
  \left( \frac{1+z_1z_2}{2} \right)^n, \quad k,\ell,n \in \NN_0\,.
$$
In the case $k=\ell=n=2$, results in \cite{CharinaContiJetterZimm2010} imply that the degree of
polynomial generation is $k=4$. Now, to check the degree of
polynomial reproduction we first use Lemma \ref{lemma:box} to compute
$\btau$ that guarantees the reproduction of linear polynomials, \ie $\
\btau=\frac{1}{2}(k+n, \ell+n)=(2,2)$. Using Theorem \ref{teo:main} with this
$\btau$ we see that the scheme does not reproduce polynomials of
degree $k=2$, since $q_{1,1}(\btau)=\btau_1 \cdot \btau_2=4$, but
$ D^{(1,1)} B_{2,2,2}(\bone)=18$.

\subsection{Three dimensional example}

Three dimensional examples can also be considered. For example, it is easy to show that  the subdivision  scheme  with the mask symbol
\begin{eqnarray*}
 a(z_1,z_2,z_3)&=&2^3\Big[6z_1z_2z_3 \left(\frac{1+z_1}{2}\right)^2\left(\frac{1+z_2}{2}\right)^2\left(\frac{1+z_3}{2}\right)^2
\left(\frac{1+z_1z_2z_3}{2}\right)^2 - \\
&& \frac{5}{4}z_1
\left(\frac{1+z_1}{2}\right)\left(\frac{1+z_2}{2}\right)^3\left(\frac{1+z_3}{2}\right)^3
\left(\frac{1+z_1z_2z_3}{2}\right)^3 -\\
&& \frac{5}{4}z_2
\left(\frac{1+z_1}{2}\right)^3\left(\frac{1+z_2}{2}\right)\left(\frac{1+z_3}{2}\right)^3
\left(\frac{1+z_1z_2z_3}{2}\right)^3 -\\
&& \frac{5}{4}z_3
\left(\frac{1+z_1}{2}\right)^3\left(\frac{1+z_2}{2}\right)^3\left(\frac{1+z_3}{2}\right)
\left(\frac{1+z_1z_2z_3}{2}\right)^3 -\\
&& \frac{5}{4}z_1
z_2z_3\left(\frac{1+z_1}{2}\right)^3\left(\frac{1+z_2}{2}\right)^3\left(\frac{1+z_3}{2}\right)^3
\left(\frac{1+z_1z_2z_3}{2}\right)\Big].
\end{eqnarray*}
is convergent. Moreover,
using Theorem \ref{teo:main} we get that for $\btau=(3,3,3)$, thus,  the scheme also reproduces linear polynomials. Since
$ D^{(2,0,0)} a(\bone)=46 \not = 8 \cdot q_{2,0,0}(\btau)=48$, the scheme is not reproducing polynomials of degree $2$.

\subsection{Interpolatory schemes}
Exactly the same argument for interpolatory schemes as in
\cite[Corollary 5.3]{ContiHormann} extends to the multivariate
case. Interpolatory schemes are such whose mask satisfies
$$
 \ra_\bnull=1, \quad \ra_{m\balpha}=0, \quad \balpha \in \ZZ^s\,,
$$
and therefore, when convergent, with a limit function interpolating the initial data as well as all the data
generated through the recursions.

\smallskip \noindent
Let us assume that an interpolatory scheme generates polynomials
up to degree $k$. Due to the special structure of the symbol of
interpolatory schemes
$$
  a(\bz)=1+\sum_{\be \in E \setminus \{\bnull\}} a_\be(\bz)
$$
we get $D^\bj a (\bone)=0$ for all $\bj \in \NN_0^s$ with
$|\bj|\le k$. Then, by Proposition \ref{prop:linear-reproduction},
the suitable choice of $\btau$ in \eqref{def:general_par} for
reproduction of linear polynomials is $\btau=\bnull$. Corollary
\ref{cor:main} and the definition of $q_{\bj}$ in \eqref{def:qj}
imply that the scheme also reproduces polynomials up to degree $k$
with this $\btau$. Thus, the following result holds and confirms that our
results reproduce results in \cite{Je}.

\begin{proposition}\label{prop:int}
 A convergent interpolatory scheme $S_\ra$ reproduces polynomials up to degree $k$ if and only
 if it generates polynomials of degree up to $k$.
\end{proposition}

For the butterfly scheme with the symbol
$$
 a(z_1,z_2)=4 \cdot z_1^{-3}z_2^{-3}[7z_1z_2B_{2,2,2}-2z_1B_{1,3,3}-2z_2B_{3,1,3}-2z_1z_2B_{3,3,1}]
$$
we have $\btau=(0,0)$, as expected. Since $a(z_1,z_2)$ satisfies sum rules of order $4$ (see, again \cite{CharinaContiJetterZimm2010}),
the subdivision scheme generates cubic polynomials, it also reproduces cubic polynomials by Corollary \ref{cor:main}.

\subsection{$\sqrt{3}-$subdivision}
The approximating $\sqrt{3}-$subdivision scheme from \cite{JiangOswald} with the mask symbol
\begin{eqnarray}
 a(\bz)&=&\frac{1}{6}\left(z_1z_2+z_1^{-1}z_2^{-1}+z_1^{-1}z_2^2+z_1^{-2}z_2+z_1z_2^{-2}+z_1^2z_2^{-1}\right)
 \notag \\
 &+&\frac{1}{3} \left(z_1^{-1}+z_2+z_1z_2^{-1}\right)+\frac{1}{3}
 \left(z_2^{-1}+z_1+z_1^{-1}z_2\right)\notag
\end{eqnarray}
satisfies sum rules at most of order $2$. The associated dilation matrix $M=\left[\begin{array}{rr} 1&2\\-2&-1 \end{array}\right]$ satisfies $M^2=-3I$ and the corresponding refinable function is also refinable
with respect to the iterated mask $a(z_1z_2^{-2},z_1^2z_2^{-1}) \cdot a(\bz)$. By Corollary \ref{cor:linear-reproduction}, the
corresponding scheme reproduces linear polynomials, if $\btau=(0,0)$. Thus, associated refinable function
has approximation order $2$.

\section{Conclusions}
In this paper we give algebraic conditions on the symbol  of a multivariate
subdivision scheme with dilation matrix $mI,\ |m|\ge 2$, that allow us to determine
the degree of polynomial reproduction of the scheme.
These conditions also yield the correct parametrization
for any convergent  subdivision scheme to guarantee polynomial
reproduction of degree at least $1$. This is true in particular
for subdivision schemes associated with box splines. The
restriction of a dilation matrix of type $mI$ and tensor product
structure of the polynomial in (\ref{def:qj}) let us extend the
univariate results in \cite{ContiHormann} easily to the multivariate
setting. We believe that this paper is an important
first step towards the investigation of polynomial reproduction of
multivariate subdivision schemes with general dilation matrix, which is
currently under investigation.

\bigskip
\bibliographystyle{amsplain}

\vspace{10mm}

\noindent\hbox{\vtop{\hsize15pc\parindent0pt
Maria Charina\\
Fakult\"at f\"ur Mathematik\\
TU Dortmund\\
Vogelpothsweg 87\\
D--44227 Dortmund\\
{\tt maria.charina@uni-dortmund.de} }}
\hskip3pc
\noindent\hbox{\vtop{\hsize15pc\parindent0pt
Costanza Conti\\
Dipartimento di Energetica\\
Universit\`a di Firenze\\
Via C. Lombroso 6/17\\
I--50134 Firenze\\
{\tt costanza.conti@unifi.it} }}

\end{document}